\documentclass[12pt,reqno]{amsart}
\usepackage{amsmath, amsfonts, amssymb, amsthm}

\newcommand{\Ff}{{\mathbb F}}
\newcommand{\Zz}{{\mathbb Z}}

\newcommand{\Qq}{{\mathbb Q}}

\newcommand{\Ll}{{\mathcal L}}

\newcommand{\Hh}{{\mathcal H}}

\newcommand\bb{{\mathcal B}}        %

\newtheorem{Theorem} {Theorem}% [section]
\newtheorem{Proposition} [Theorem] {Proposition}
\newtheorem{Lemma} [Theorem] {Lemma}
\newtheorem{Corollary} [Theorem] {Corollary}

\newtheorem{Definition}[Theorem]{Definition}

\newcommand{\Proof}{ \noindent{\bf Proof:}\quad }

\def\Tr{\operatorname{Tr}}

\def\PG{\operatorname{PG}}
\def\AG{\operatorname{AG}}

\def\eqref#1{(\ref{#1})}

\begin{document}

\title[unital intersection sizes]
{
  The sizes of the intersections of two unitals in $\PG(2,q^2)$}
\author{David B. Chandler${}^*$}

\keywords{unital, Hermitian variety, incidence matrix, $p$-adic ring}
\thanks{${}^*$This work was partially supported by NSF grants  DMS 
0701049 and DMS 1001557.}
%\address{Institute of Mathematics, Academia Sinica,
%Nankang, Taipei  11529,  Taiwan, Republic of China}
\email{davidbchandler@gmail.com}
 %\address{Department of Mathematical Sciences, University of 
%Delaware,\\
 % Newark, DE 19716, USA

%\keywords{}\subjclass{ 05E20 (Primary), 20G05 20C11 (Secondary)}
%\renewcommand{\subjclassname}{\textup{2000} Mathematics Subject 
%\Classification}
\begin{abstract}
We show that the size of the intersection of a Hermitian variety 
in $\PG(n,q^2)$, and any set satisfying an $r$-dimensional-subspace
intersection property, is congruent to 1 modulo a power of $p$.  
In particular, in the case where $n=2$, if the two sets are a Hermitian 
unital and any 
other unital, the size of the intersection  is congruent to 1 modulo $\sqrt 
q$ or modulo 
$\sqrt{pq}$.  If the second unital is a Buekenhout-Metz unital, 
we show that the size is congruent to 1 modulo $q$.

\end{abstract}

\maketitle
\section{Introduction}\label{1}
A \emph{unital} is any 2-design with parameters of the form 
$(n^3+1,n+1,1)$.  That is, we have a set $\mathcal P$ of $v=n^3+1$ points 
and a collection $\mathcal B$, the blocks, of subsets of $\mathcal P$, 
having the following two properties:  Each block has size $k=n+1$; and any 
two points are jointly contained in exactly one block.  A unital $U$ is 
\emph{embedded} in a finite projective plane of order $q^2$ if it is a set 
of $q^3+1$ points of the plane with the property that every line of the 
plane intersects exactly 1 or $q+1$ points of $U$.  In this paper we are 
interested in unitals which are embedded in 
$\PG(2,q^2),\ q=p^t$, for some prime $p$.

The classical or \emph{Hermitian} unital is a Hermitian variety, the set of 
zeroes of a unitary form.  Such forms on
$\Ff^3_{q^2}$ are projectively equivalent to
\begin{equation}%\label{Hequation}
x^{q+1}+y^{q+1}+z^{q+1}
\nonumber\end{equation}
over $\Ff_{q^2}$, the field with $q^2$ elements.  More generally, we are also interested in 
 Hermitian varieties in $\PG(n,q^2)$, which are 
projectively equivalent to $\sum_{i=0}^n x_i^{q+1}=0$.

Buekenhout and Metz proved the existence of nonclassical unitals in the  
1970s \cite{bm}.  
Around 1990 Baker and Ebert generalized Buekenhout's and Metz's construction  to describe a 
two-parameter family of unitals $U_{a,b}=\{(0,1,0)\}\cup\{(x,ax^2+bx^{q+1}+r,1):
x\in\Ff_{q^2}\ \mathrm{and}\ r\in\Ff_q\}$, in any 
desarguesian plane of 
square order, where $a,b\in\Ff_{q^2}$ meet some condition and $q>2$ \cite{numeratium,bakerebert,gary}.    We call these unitals \emph{Buekenhout-Metz unitals}, or \emph{B-M unitals}.  A B-M unital is Hermitian if and only if $a=0$.  The construction is summarized in Section~\ref{bmun}.

  When $q$ is is an odd power of 2, there is one nonclassical ovoid known in 
$\PG(3,q)$, 
the \emph{Tits ovoid}.  By 
replacing the elliptic quadric in the B-M construction by a Tits ovoid, 
we get one more projective equivalence class of unitals in $\PG(2,q^2)$.  
At present there are no other unitals known in desarguesian planes.  For 
more information about unitals, we refer the reader to
\cite{be}.

Kestenband \cite{kesten} showed that if $H_1$ and $H_2$ are Hermitian 
unitals in $\PG(2,q^2)$, then $|H_1\cap 
H_2|\in\{1,q+1,q^2-q+1,q^2+1,q^2+q+1,q^2+2q+1\}$, and also determined the 
possible intersection configurations.  Note that the sizes are all 
congruent to 1 modulo $q$.  Baker and Ebert  \cite{numeratium} then proved 
that the size of 
the intersection of a Hermitian unital with the special type of B-M unital 
having $b=0$ is congruent to 1 modulo $p$.  They also 
conjectured that the size of the intersection of any unital with a 
Hermitian unital would turn out to be congruent to 1 modulo $q$.  Blokhuis, 
Brouwer, and Wilbrink \cite{bbw} soon proved that the size of this 
intersection is 
congruent to 1 modulo $p$, by showing that the Hermitian unital is in the 
code of lines of $\PG(2,q^2)$.    In this article we prove the 
following theorem.
\begin{Theorem}\label{unitalint}
Let $H$ be a Hermitian unital embedded in $\PG(2,q^2)$ and let $U$ be any 
other unital embedded in $\PG(2,q^2)$, where $q=p^t$, and $p$ is a prime.  
Then
$$|H\cap U|\equiv1\quad (\mathrm{mod}\ p^{\lceil t/2\rceil}).$$
Moreover, if $U$ is a Buekenhout-Metz unital, then
$$|H\cap U|\equiv1\quad (\mathrm{mod}\ q).$$
\end{Theorem}

We should note that if neither $U_1$ nor $U_2$ is Hermitian, then nothing 
particular can be said about $|U_1\cap U_2|$.  A computer check of the intersection 
sizes of pairs randomly chosen from the known non-Hermitian unitals reveals no particular pattern.

\section{Hermitian \textit{vs.} arbitrary unital}\label{HvA}
In this section we obtain the first part of Theorem~\ref{unitalint} as a corollary of a more general result.

Let $p$ be a prime, let $q=p^t$, and let $V$ be an $(n+1)$-dimensional 
vector space  
over $\Ff_{q^2}$ with coordinate functions $x_0,\ldots,x_n$.  We denote the 
set of 
projective points of $\PG(n,q^2)$, {\it i.e.}, one-dimensional subspaces of $V$, by 
$\Ll_1$, 
and the set of projective $(r-1)$-dimensional subspaces $\PG(n,q^2)$,  {\it i.e.},
$r$-dimensional subspaces of $V$, by $\Ll_r$.  

Let $H$ be a Hermitian variety of $\PG(n,q^2)$.  Note that every Hermitian variety is 
projectively equivalent to the zeroes of
$$\sum_{i=0}^n x_i^{q+1}.$$
Suppose we have $r$, a
vector-space dimension, $1<r\le n$,  another  positive integer $\beta$,
and a set of points $S \subset \Ll_1$ with the following 
intersection property:
\\

\noindent {\bf Property I.}
  Every element of $\Ll_r$ meets $S$ in 
a number of points which is divisible by $p^\beta$.
\\

\noindent Then we will prove that $|S\cap H|$ is divisible by a certain 
power of $p$.  Note that from the size of a projective $(r-1)$-space,
$$\beta\le 2t(r-1) \le 2t(n-1).$$
\begin{Definition}
Let $A_{r,1}$ be the $(0,1)$-matrix, columns  indexed by $\Ll_1$, the points, and rows by $\Ll_r$, the projective $(r-1)$-spaces, whose entries are
$$a_{Y,Z}=\left\{\begin{array}{cl}
1, & \mathrm{if}\ Z\subset Y;\\
0, & \mathrm{otherwise};
\end{array}\right.\quad Y\in\Ll_r,\ Z\in\Ll_1.$$
We call $A_{r,1}$ the \emph{incidence matrix} between $\Ll_r$ and $\Ll_1$.
\end{Definition}

Consider $\Ff_{q^2}^{\Ll_1}$, the space of $\Ff_{q^2}$-valued functions on $\Ll_1$.
A useful  basis for this space is given by the set of monomials \cite{bsin}:
$$\overline{\mathcal B}=\left\{\prod_{i=0}^n
x_i^{b_i},\quad
{\begin{array}{c}
0\le b_i\le q^2-1,\ 0\le i\le n,\quad
(q^2-1) \mid \sum_{i=0}^n b_i,\\
(b_0,\ldots,b_n)\ne (q^2-1,\ldots,q^2-1)
\end{array}}\right\}.$$

It
will be helpful to view the entries of
 $A_{r,1}$ as coming from some $p$-adic local ring. Let $q=p^t$ and let
$K=\Qq_p(\xi_{q-1})$ be the unique unramified
 extension of degree $t$ over $\Qq_p$, the field of $p$-adic numbers,
where $\xi_{q-1}$ is a primitive $(q-1)^{\rm
 th}$ root of unity in $K$. Let $R=\Zz_p[\xi_{q-1}]$ be the ring of
integers in $K$ and let $\mathfrak p$ be the
 unique maximal ideal in $R$. Then 
the reduction of $R \
 (\mathrm{mod}\,\mathfrak p)$ will be $\Ff_q$. Define $\bar x$ to be $x\
(\mathrm{mod}\,\mathfrak p)$ for $x\in
 R$. Let $T_q$ be the set of roots of $x^q=x$ in $R$ (a Teichm\"uller set)
and let $T$ be the Teichm\"uller
 character of $\Ff_q$, so that $T(\bar x)=x$ for $x\in T_q$. 
We adopt the convention that $T^0(\bar x)=1,\ \bar x\in\Ff_q $, while
$T^{q-1}(0)=0$, and $T^{q-1}(\bar x)=1, \ \bar x\in\Ff_q^*$.
We will use
$T$ to lift a basis of $\Ff_q^{\mathcal
L_1}$ to a basis of $R^{\mathcal L_1}$.

For any $u\in 
R$, we 
define $\nu_p(u)$ 
to be the $p$-adic valuation of $u$.  That is, $\nu_p(u)=\alpha$ if 
$p^\alpha\mid u$ but 
$p^{\alpha+1}\nmid u$.

We obtain  a \emph{lifted}  
basis $\bb$ for the free module $R^{\Ll_1}$ (see \cite{csx} for proof):
$$\bb=\{T(\prod x_i^{b_i})\mid \prod x_i^{b_i}\in \overline\bb\}.$$

The matrix $A_{r,1}$ can be viewed as a map from 
$\Ff_{q^2}^{\Ll_1}$ to $\Ff_{q^2}^{\Ll_r}$
(or \textit{vice versa}).  
For instance, let $\mathbf{u}$ be the column (0,1)-characteristic vector of 
a point set.  Then $A_{r,1}\mathbf{u}$ records the number of points in each 
$(r-1)$-dimensional subspace of $\PG(n,q^2)$.  

In \cite{csx} it was shown that $\bb$ forms a ``Smith normal form" basis 
for the map from $\Ll_1$ to $\Ll_r$:  Let $\mathbf{v}$ be the column 
vector representing an element of $\bb$.  Then $A_{r,1}\mathbf v$
 is the corresponding invariant (a power of $p$) multiplied by an integral vector 
indexed by $\Ll_r$.

We recall the formula for the 
invariants (stated for $\Ff_{q^2})$.
  To each nonconstant basis monomial 
$f=x_0^{b_0}\cdots x_n^{b_n}\in\overline\bb$,
 we associate a pair of $2t$-tuples, 
$(\lambda_0,\ldots,\lambda_{2t-1})$ (called the \emph{type}), and $(s_0,\ldots,s_{2t-1})$
(called the \emph{$\Hh$-type}).  The type of $T(f)$ is that of $f$.  We expand each exponent as
$$b_i=a_{i,0}+a_{i,1}p+\cdots+a_{i,2t-1}p^{2t-1};\quad
\begin{array}{l} 0\le a_{i,j}\le p-1,\\
0\le i\le n,\\ 0\le j\le 2t-1.
\end{array}$$
Then we define
\begin{equation}\label{lambda}
\lambda_j=a_{0,j}+a_{1,j}+\cdots+a_{n,j}
\end{equation}
\begin{eqnarray}\nonumber
s_j&=&\frac1{q^2-1}\sum_{i=0}^n\left(\sum_{\ell=0}^{j-1}
p^{\ell+2t-j}a_{i,\ell}+\sum_{\ell=j}^{2t-1}p^{\ell-j}a_{i,\ell}
\right)\\ \label{ssum}
&=&\frac1{q^2-1}\sum_{i=0}^n\left(p^{2t-j}b_i\ (\mathrm{mod}\ {q^2-1}) \right)
\end{eqnarray}
and we have the relation $\lambda_j=ps_{j+1}-s_j$ (subscripts modulo $2t$).  The numbers $(q^2-1)s_j$ are called the \emph{twisted degrees} of $f$.  The formula we want is given as follows.
\begin{Proposition}[\cite{csx}]\label{invariants}
Let $f\in\bb$ be a basis monomial.  If $f=1$, the corresponding $p$-adic 
invariant is 1.  Otherwise let the
$\Hh$-type of $f$ be $(s_0,\ldots,s_{2t-1})$.  Then the  corresponding 
$p$-adic  invariant for 
the map $A_{r,1}$ from $\Ff_{q^2}^{\Ll_1}$ to $\Ff_{q^2}^{\Ll_r}$ is given 
by 
$p^\alpha$, where
$$\alpha=\sum_{j=0}^{2t-1}\max\{0,r-s_j\}.$$
\end{Proposition}

For any integer $u$, we also define define $\sigma_p(u)$ to be the  
$p$-adic digit sum.  That is, if
$$u=\sum_{j=0}^m u_j p^j,\ 0\le u_j\le p-1,\ 0\le j\le m,\quad\mathrm{then}
\quad \sigma_p(u)=\sum_{j=0}^m u_j.$$
 Note that 
$\sigma_p(u_1u_2)\le\sigma_p(u_1)\sigma_p(u_2),\ 
\sigma_p(u_1+u_2)\le\sigma_p(u_1)+\sigma_p(u_2)$,
and
\begin{equation}\label{digitsums}
\sum_{i=0}^n 
\sigma_p(b_i)=\sum_{j=0}^{2t-1}\lambda_j=(p-1)\sum_{j=1}^{2t-1}s_j
\end{equation}
for $f=x_0^{b_0}\cdots x_n^{b_n}$ of type $(\lambda_0,\ldots,\lambda_{2t-1})$.  
Here, $b_i$ is taken as the least positive residue modulo $(q^2-1)$, unless 
it is already zero.  Note that reduction to the least positive residue may 
reduce $\sigma_p(b_i)$, but never increase it.

 We use the fact that
\begin{equation}\nonumber%\label{galois}
T(a+b)\equiv \left(T(a)+T(b)\right)^{q^{\ell}}\quad (\mathrm{mod}\ 
q^{\ell}),\quad\mathrm{for}\ a,b
\in\Ff_{q},
\end{equation}
for any positive integer $\ell$, which enables us to bring the Teichm\"uller 
character inside the parentheses \cite{csx}.

Let $H$ be the Hermitian variety given by
$$\sum_{i=0}^n x_i^{q+1}=0.$$
Then we take $\overline H$ as the complement of $H$ in 
$\PG(n,q^2)$.  Consider 
the characteristic function:
\begin{eqnarray}\nonumber
\chi_{\overline H}:\ \Ll_1&\to&\{0,1\}\subset R\\
\nonumber\chi_{\overline 
H}\left(\langle(x_0,x_1,\ldots,x_n)\rangle\right)&=&T\left(\sum_{i=0}^n 
x_i^{q+1}\right)^{q-1}\\
&\equiv& \left(\sum_{i=0}^n 
T(x_i)^{q+1}\right)^{q^{2\ell+1}-q^{2\ell}}\label{Hfunc}
\quad\mathrm{mod}\ q^{2\ell}.
\end{eqnarray} 
We see that each term in the expansion of 
the RHS of (\ref{Hfunc})
has the form
\begin{equation}\nonumber %\label{Hterm}
f={q^{2\ell+1}-q^{2\ell}\choose k_0,\ldots,k_n}
\prod_{i=0}^n T(x_i)^{(q+1)k_i}
\end{equation}
where $\sum_{i=0}^n k_i=q^{2\ell+1}-q^{2\ell}$.  From Legendre's formula, 
$\nu_p(n!)=\frac{n-\sigma_p(n)}{p-1}$, we have
$$(p-1)\/\nu_p{q^{2\ell+1}-q^{2\ell}\choose k_0,\ldots,k_n}=\sum_{i=0}^n 
\sigma_p(k_i)
-(p-1)t.$$
Since $$2\sigma_p(k_i) \ge \sigma_p\left((q+1)k_i\right)\ge \sigma_p(b_i), $$ 
where $b_i$ is the least positive residue of $(q+1)k_i$ modulo $(q^2-1)$, if 
$k_i>0$, we substitute from (\ref{digitsums}) to obtain 
\begin{equation}\label{hermsum}
2\,\nu_p{q^{2\ell+1}-q^{2\ell}\choose k_0,\ldots,k_n}\ge \frac1{p-1}
\sum_{j=0}^{2t-1}\lambda_j-2t=
\sum_{j=0}^{2t-1}s_j-2t
\end{equation} 
where $(\lambda_0,\ldots,\lambda_{2t-1})$ and $(s_0,\ldots,
s_{2t-1})$ are the type and $\Hh$-type of $f$.  In the special case  where
$b_0=b_1=\cdots=b_n=(q^2-1)$, the $\Hh$-type is actually not defined in 
\cite{bsin}, but the same calculation shows that $q^{n}$ divides the 
multinomial coefficient in this case (as if the $\Hh$-type were 
$(n+1,\ldots,n+1)$).

We now consider the set $S$ with Property I.
Let 
\begin{eqnarray*}%\label{Sfunc}
\chi_S:\ \Ll_1&\to&\{0,1\}\subset R \\
\chi_S(\langle(x_0,x_1,\ldots,x_n)\rangle)&=&
\chi_S(x_0,x_1,\ldots,x_n)\\
\chi_S(\langle(0,\ldots,0)\rangle)&=&\chi_S(0,\ldots,0)=0
\end{eqnarray*}
be the characteristic function of 
$S$ expressed as a polynomial in $T(x_0)$, $\ldots$, $T(x_n)$.  
We now assume that $\chi_S$ is restricted to points of $V$ other than the origin 
and use the identity,
\begin{equation}\label{offorigin} 
\prod_{i=0}^n (1-x_i^{q^2-1})=0,
\end{equation}
to eliminate the monomial $\prod_{i=0}^n x_i^{q^2-1}$.
Let 
$g=c_S\prod_{i=0}^n T(x_i)^{b'_i}$ be a monomial term of 
$\chi_S$, where $c_S\in R$, and the
$\Hh$-type of $g$ be $(s'_0,\ldots,s'_{2t-1})$. 
Since the monomials form a Smith-normal-form basis for the incidence 
matrix $A_{r,1}$, our divisibility property (\emph{i.e.}, Property I) implies that 
$p^\beta$ must 
divide the product of $c_S$ and the invariant corresponding to $g$ (for 
the matrix $A_{r,1}$). 
Thus $p^\beta$ divides the constant term, and for each nonconstant term $g$, we 
get
$$\nu_p(c_S) + \sum_{j=0}^{2t-1}\max\{0,r-s'_j\}\ge\beta,$$
or
\begin{equation}\label{usum}
\nu_p(c_S) \ge 
\max\left\{0,\beta-2rt+\sum_{j=0}^{2t-1}\min\{r,s'_j\}\right\}.
\end{equation}

Observe that 
$$|S\cap \overline H|=\frac 1 {q^2-1} 
\sum_{\mathbf x=(x_0,x_1,\ldots,x_n)\in V} \chi_S({\bf x})\chi_{\overline 
H}({\bf x})\,,$$
since the function we have for $\chi_{\overline H}$ evaluates to zero at 
the origin.  Let 
$$f=c_{\overline H}\prod_{i=0}^n T(x_i)^{b_i}, \quad 
c_{\overline H}={q^{\ell+1}-q^\ell \choose k_0,\ldots, k_n},$$ 
be some term in the expansion of the RHS of (\ref{Hfunc}) with $\Hh$-type\\ 
$(s_0,\ldots,s_{2t-1})$, and let 
$$g=c_S\prod_{i=0}^n T(x_i)^{b'_i}$$ be some monomial term of $\chi_S$ with 
$\Hh$-type $(s'_0,\ldots,s'_{2t-1})$.  We wish to show that a certain 
power of $p$ always divides
\begin{eqnarray*}
\sum_{\mathbf x\in V}fg&=&c_{\overline H}c_S\sum_{\mathbf x\in V}\,
\prod_{i=0}^n T(x_i)^{b_i+b'_i}\\
&=&c_{\overline H}c_S \prod_{i=0}^n \sum_{x_i\in\Ff_{q^2}}
T(x_i)^{b_i+b'_i}.
\end{eqnarray*}

Since
\begin{equation}
\label{monomialsum}
\sum_{x\in\Ff_{q^2}}T(x)^j=\left\{
\begin{array}{ccl}
q^2,&\ \mathrm{if}\ & j=0,\\
0,&\ \mathrm{if}\ & j \not\equiv 0\ (q^2-1),\\
q^2-1,&\ \mathrm{if}\ & j>0,\ j \equiv 0\ (q^2-1),
\end{array}
\right.
\end{equation} 
we only need to consider terms with $b_i+b'_i\in \{0,q^2-1,2(q^2-1)\},\ 
0\le i \le n$. 
We first suppose that $b_i+b'_i=q^2-1$ for $0\le i\le n$.  In this case, 
$fg$ has total degree 
$(q^2-1)(n+1)$, and all the twisted degrees also are $(q^2-1)(n+1)$.  
Therefore $s_j+s'_j= n+1,\ 
0\le j\le 2t-1$.  Now,

\begin{eqnarray}
\label{rhs} 
2\nu_p(c_{\overline H}c_S) & \ge & \sum_{j=0}^{2t-1} (s_j-1)+
2\max\left\{0,\beta-2rt+\sum_{j=0}^{2t-1}\min\{r,s'_j\}\right\}\ {}\\&
= & \sum_{j=0}^{2t-1} (n-s'_j)+
\max\left\{0,2\beta-4rt+2\sum_{j=0}^{2t-1}\min\{r,s'_j\}\right\}.
\nonumber\end{eqnarray}

Notice that increasing one $s'_j$ by one decreases the sum by one if $s'_j 
\ge r$ or if the second term inside the ``$\max$'' function of (\ref{rhs})
 is negative, and it 
increases the sum by one otherwise.  Therefore the minimum value of the 
sum can be achieved by choosing all the $s'_j$ to be either 1 or $n$, 
except for one value, say $s'_0$, to be in the range $1<s'_0<r$.  If
$r-s'_j\le n-r$, the sum does not increase if we change the value to 
$s'_j=n$.
Therefore,
if $r\le (n+1)/2$, the right-hand 
side is minimized if we always pick $s'_j=n$. 
If $r>(n+1)/2$, the 
right-hand side is minimized in one of two ways: either we choose 
$s'_j\in \{1,n\},\ j>0$, and $s'_0$ as necessary to make the second term 
on 
the RHS equal to zero, or if $s'_0\ge 2r-n$ in that case, we also make 
$s'_0=n$.

We also consider the cases where $b_i+b'_i\in \{0,2(q^2-1)\}$ for some $i$.  
If $b_i=b'_i=0$, then $s_j+s'_j$ is reduced by 1 for each $j,\ 0\le j \le 
2t-1$, compared to the case where $b_i=q^2-1$ and everything else is the 
same, 
which reduces our estimate for $2\nu_p(c_{\overline H}c_S)$ by $2t$, 
but from (\ref{monomialsum}), we have an extra factor of $q^2$, which 
increases our estimate by $4t$.  If $b_i=b'_i=q^2-1$, then each $s_j $ is 
increased by 1, compared to the case where $b_i=0$, which only increases our 
estimate.

Summarizing the results for  $r>(n+1)/2$, we get the smallest estimate by 
choosing the number of $j$ such 
that $s'_j=1$ to be
 $\alpha=\left\lfloor \frac {\beta} {r-1} \right\rfloor$ and
let
$\gamma=\beta - (r-1)\alpha$.  Then $s'_0\in \{r-\gamma,n\}$.

  Since Property I
 also implies $\nu_p(|S|)\ge \beta$, we have proved the following:

\begin{Theorem}
Let $S\subset\Ll_1$ be a set of points with Property I,
and let $H\subset \Ll_1$ be the point set of a nondegenerate 
Hermitian variety.  Let $\alpha$ and $\gamma$ be as above and let
$$\theta=\left\{\begin{array}{lcc}
\beta,&\mathrm{if}& r\le (n+1)/2; \\
\left\lceil \frac{(n-1)\alpha} 2 + \min\left\{\frac
{n-r+\gamma}
2,\gamma\right\} \right\rceil, & \mathrm{if} & r> (n+1)/2 .
\end{array}\right.$$
Then $$p^\theta \mid |S\cap H|.$$
\end{Theorem}
We note that $\theta$ is approximately $\frac{\beta(n-1)}{2(r-1)}$ if
$r> (n+1)/2$.

Taking $n=r=2,\ \beta=t$, and $S=\overline U$, we have:
\begin{Corollary}
Let $H$ be a Hermitian unital and let $U$ be an arbitrary unital in $\PG(2,q^2),\ 
q=p^t$.  
Then $\nu_p(|H\cap U|-1)\ge t/2$.
\end{Corollary}

%%%%%%%%%%%%%%%%%%%%%%%%%%%%%%%%%%%%
\section{Hermitian \textit{vs.} Buekenhout-Metz unital}\label{bmun} In this section we show that 
the number of points in the intersection of a Hermitian  unital 
and a Buekenhout-Metz unital is always congruent to 1 (modulo $q$).
 The Buekenhout-Metz construction goes as follows (see 
\cite{bm,gary,bakerebert}).  Start with an 
elliptic quadric in $A\cong\PG(3,q)\subset\PG(4,q)$ and a 
regular spread in 
$B\cong\PG(3,q)\subset\PG(4,q)$ such that the plane $A\cap B$ 
is a tangent plane to the 
quadric.  Choose a point $P$ on the same spread line as the 
point of tangency of $A\cap B$ 
to the quadric and let $U^*$ be the cone of $P$ and the 
quadric.  Now construct a new 
plane of order $q^2$, taking as points the points of 
$\PG(4,q)\setminus B$ as well as the 
spread lines covering $B$, and taking as lines the planes of 
$\PG(4,q)$ whose intersection 
with $B$ is a spread line, as well as $B$ itself.  With 
inclusion as incidence, we always 
get a translation plane 
(having a regular automorphism group on the image of 
$\PG(4,q)\setminus B$), and since we 
took a regular spread we get the desarguesian plane 
$\PG(2,q^2)$.  Furthermore, the 
cone $U^*$ contains $q^3$ affine points and a spread line. The 
image $U\subset\PG(2,q^2)$ 
of $U^*$ is easily shown to be a unital.

Every unital produced by the Buekenhout-Metz construction is 
projectively equivalent to
$$U_{\alpha,\beta}=\left\{(1,y,\alpha y^2+\beta y^{q+1}+r)\mid 
y\in\Ff_{q^2},r\in\Ff_q\right\}\cup\{(0,0,1)\}.$$
If $\alpha=0$ and $\beta\notin\Ff_q$ then $U_{0,\beta}$ is a 
Hermitian 
unital.  If $\alpha\neq0$ then 
$U_{\alpha,\beta}$ is a unital if and only if the following 
condition holds:
\begin{eqnarray*}
(\beta^q-\beta)^2+4\alpha^{q+1}\textrm{ is a nonsquare of } 
\Ff_q,\quad&\mathrm{if}&q\textrm{ is odd;}\\
\Tr_{\Ff_q/\Ff_2}\left(\frac{\beta^q+\beta}{\alpha^{q+1}} 
\right)=1,\quad\quad&\mathrm{if}&q\textrm{ is 
even.}\\
\end{eqnarray*}
From $z-(\alpha y^2+\beta y^{q+1})=r\in\Ff_q$ (if $x=1$) we get 
the 
equation
\begin{equation}\label{bmeq}
\alpha^qy^{2q}-\alpha y^2+(\beta^q-\beta)y^{q+1}-z^q+z=r-r^q=0
\end{equation}
which is satisfied by the affine points of $U_{\alpha,\beta}$. 
The homogeneous equation for the affine points is
$$\alpha^qx^2y^{2q}-\alpha 
x^{2q}y^2+(\beta^q-\beta)x^{q+1}y^{q+1}
-x^{q+2}z^q+x^{2q+1}z=0.$$
 Note 
that if $q$ is even, then the 
left-hand side of (\ref{bmeq}) is in $\Ff_q$, while if $q$ is 
odd, then the square of the 
left-hand side of (\ref{bmeq}) is in $\Ff_q$.  Therefore in 
either case, raising to the $2(q-1)$ power, we 
get the characteristic function over $\Ff_{q^2}$ for the affine 
points of the complement of $U_{\alpha,\beta}$.
The complete characteristic function of the complement of 
$U_{\alpha,\beta}$, including all the points at infinity 
other than 
$(0,0,1)$, is
\begin{eqnarray*}
&\left(\alpha^qx^2y^{2q}-\alpha 
x^{2q}y^2+(\beta^q-\beta)x^{q+1}y^{q+1}
-x^{q+2}z^q+x^{2q+1}z\right)^{2(q-1)}
&\\&+\,(1-x^{q^2-1})y^{q^2-1}.&
\end{eqnarray*}

   Taking the $p$-adic Teichm\"uller 
character (modulo $q^2$) of each term inside the parentheses
and further raising the affine part of the characteristic 
function  to the $q^2$ power, we finally   get the 
characteristic function $\chi_{\overline U}$ in the Galois ring 
$R/q^2R$ ({\it cf.} \cite{wan,csx}).

We also need the characteristic function of the Hermitian 
unital.  Previously the 
orientation of the non-Hermitian unital was arbitrary; so we 
could take the simplest form 
of the equation for the Hermitian unital.  Now that the 
orientation of the B.-M. unital is 
fixed, we must consider a general form for the Hermitian unital:
\begin{eqnarray}
(c_{1,1}x+c_{1,2}y+c_{1,3}z)^{q+1}+(c_{2,1}x+c_{2,2}y+c_{2,3}z)^{q+1}
&&\nonumber\\
\nonumber+\;(c_{3,1}x+c_{3,2}y+c_{3,3}z)^{q+1}&=&\\
(c_{1,1}^qx^q+c_{1,2}^qy^q+c_{1,3}^qz^q) \label{generalHerm}
\cdot(c_{1,1}x+c_{1,2}y+c_{1,3}z)&&\\ \nonumber
+\ (c_{2,1}^qx^q+c_{2,2}^qy^q+c_{2,3}^qz^q) 
\cdot(c_{2,1}x+c_{2,2}y+c_{2,3}z)&&\\
+\ (c_{3,1}^qx^q+c_{3,2}^qy^q+c_{3,3}^qz^q) \nonumber
\cdot(c_{3,1}x+c_{3,2}y+c_{3,3}z)&=&0,
\end{eqnarray}
where the coefficients $(\mathbf c)_{(3\times3)}$ form a 
nonsingular matrix over $\Ff_{q^2}$.
Again we take the Teichm\"uller character of $x$, $y$, $z$, and 
$(\mathbf c)$ 
and raise to the $(q^3-q^2)$ 
power (since $q+1$ powers of $\Ff_{q^2}$ elements are already in 
$\Ff_q$)
 to get  $\chi_{\overline H}$, the characteristic function of 
the complement of 
the Hermitian unital in the ring
$R/q^2R$.
 
Our goal is to show that $q$ divides
$$\sum_{(x,y,z)\in{\Ff_{q^2}}^3}\chi_{\overline H}(x,y,z)\;\chi_{\overline 
U}(x,y,z)=(q^2-1)\ 
|\overline H\cap\overline U| .$$
 As before, we do so by considering the expansion of the product.
In view of (\ref{monomialsum}),
we only need to consider terms 
which arise as the product of a monomial term ${
f=c_HT(x^{b_{1}}y^{b_{2}}z^{b_{3}})}$ from the expansion of $\chi_{\overline 
H}$
and a term $g=c_UT(x^{b'_{1}}\protect{}
y^{b'_{2}}z^{b'_{3}})$ from the expansion of $\chi_{\overline U}$,
satisfying  $b_{i}+b'_{i}\in\{0,q^2-1,2(q^2-1)\},
\ i=1,2,3$. 

If $b_i=b'_i=0$ for some $i$, we already get a factor of $q^2$ from 
(\ref{monomialsum}).  Otherwise, if $b_{i}=b'_{i}=(q^2-1)$ for some $i$, 
then the twisted degrees sum to at least $4(q^2-1)$, and either one of 
$f$ or $g$ is $T(xyz)^{q^2-1}$, or both $f$ and $g$ are of $\Hh$-type
$(2,\ldots,2)$ (if, for instance, $b_1=b'_1=q^2-1,\ b_2=b'_3,\ b_3=b'_2,$ 
and $b_2+b_3=q^2-1)$.
 In the 
first case, the coefficient of $T(xyz)^{q^2-1}$ is divisible by $q$ (eliminating 
$T(xyz)^{q^2-1}$ using (\ref{offorigin})
 makes this coefficient the 
constant term).
In the second case, both 
coefficients are divisible by $q$ (see (\ref{usum}) with $\beta=t,\ r=2,$ and 
$s_j=2,\ 0\le j<2t$).
Thus we assume $b_{i}+b'_{i}=q^2-1$.

 \begin{Lemma}
Let $(s_{0},\ldots,s_{2t-1})$ and  
$(s'_{0},\ldots,s'_{2t-1})$ be the $\Hh$-types of the two monomials $f$ and 
$g$  decribed above, such that $b_{i}+b'_{i}=q^2-1,\ i=1,2,3$.  Then 
$$s_{j}+s_{j}=3,\quad
0\le j\le 2t-1.$$
\end{Lemma}
\Proof
The corresponding twisted degrees of $f$ and $g$ always sum to $3(q^2-1)$, 
because the degree of $fg$ is invariant under Frobenius twisting.
\qed

Unlike the situation in Section~\ref{HvA}, we consider the relation between the 
coefficients $c_H$ and $c_U$ and the shapes of the $\Hh$-types $f$ and $g$.  Note 
that the $\Hh$-type is a 
$2t$-tuple consisting of 1's and 2's.

\begin{Lemma}\label{hval}
Let the $\Hh$-type of $f_H$ be $(s_0,\ldots,s_{2t-1})$ and let 
$$\xi=|\{j\mid 0\le j\le t-1,\ s_j=1, \ \textrm{and}\ 
s_{j+t}=1\}|.$$  Then $c_H$ is divisible by $p^{t-\xi}.$
\end{Lemma}
\Proof
We expand 
(\ref{generalHerm}) to get 
9 terms, each of the form $\upsilon_i=c\mu \phi^q,\ i\in\{1,\ldots,9\}$, 
where $c\in\Ff_{q^2}$, and each of 
$\mu$ and $\phi$ is one of $x$, $y$, or $z$.
Then we raise to the power $q^3-q^2$ 
and take the Teichm\"uller lifting to get
$$\chi_{\overline H}\equiv\sum_{k_1+\cdots+k_{9}=q^3-q^2}
\left(\begin{array}{c} q^3-q^2\\k_1,\ldots,k_{9}\end{array}
\right) \prod_{i=1}^{9}T(\upsilon_i)^{k_i}\ (\mathrm{mod}\ q^2).$$
From Legendre's formula we have the $p$-adic valuation of the 
multinomial coefficient is 
$$\frac1{p-1}(\sigma_p(k_1)+\cdots+\sigma_p(k_{9})- \sigma_p(q^3-q^2))= 
\frac1{p-1}(\sigma_p(k_1)+\cdots+\sigma_p(k_{9}))-t$$ where $\sigma_p$ 
again indicates 
the $p$-adic digit sum.

We consider how the digits of $k_i$ contribute to the 
$\lambda$-sums (\ref{lambda}) of the monomial. 
Let $k_i=k_{i,0}+k_{i,1}p+\cdots+k_{i,2t-1}p^{2t-1}$ and $\upsilon_i=c\mu \phi^q$.
We can think of the digit $k_{i,j}$ as contributing once to 
$\lambda_j$ (via $\mu$) and once to $\lambda_{j+t}$ (via 
$\phi$, where the subscript is 
modulo $2t$).  In fact, if there is no carry when we collect the 
exponents of $x,$ $y,$ and $z,$ with $k_{i,j}=0,\ i=1,2,3,$ whenever $j<2t$, then 
\begin{eqnarray*}
\lambda_j=\lambda_{j+t}=\sum_{i=1}^{9}(k_{i,j}+k_{i,j+t}),&& 0\le j<t,\\
s_j=s_{j+t},&& 0\le j<t,\\
\sigma(k_1)+\cdots+\sigma(k_{9}) 
&=&\lambda_0+\cdots+\lambda_{t-1} \\
&=&\lambda_t+\cdots+\lambda_{2t+1}.
\end{eqnarray*}
In this case, the lemma is a special case of (\ref{hermsum}) with equality 
throughout.

If there is a carry when we collect the exponents and reduce 
(mod $q^2-1$), say a carry from the $(j-1)^\mathrm{th}$ place to 
the 
$j^\mathrm{th}$ place, then $\lambda_{j-1}$ is reduced by $p$ 
and 
$\lambda_j$ is increased by 1, which means $s_j$ is decreased by 
1 and all the other $s$ in the type are not affected.  Since a 
carry  can only increase the
value of $\xi$ and not decrease it, the lemma holds in this case 
too.  \qed

We have a complementary lemma for the Buekenhout-Metz case.

\begin{Lemma}\label{bmval}
 Let $g=c_UT(x^{b'_{1}}y^{b'_{2}}z^{b'_{3}})$ be a term of 
$\chi_{\overline U}$ and let the 
$\Hh$-type of 
$g_U$ be 
$(s'_0,\ldots,s'_{2t-1})$ and let
$$\xi=|\{j\mid 0\le j\le t-1,\ s_j=2, \ \textrm{and}\       
s_{j+t}=2\}|.$$  Then $c_U$ is divisible by $p^{\xi}.$
\end{Lemma}
\Proof
We first consider the case that one of the exponents $b'_{1},b'_{2},$ or 
$b'_{3}$ is 0 or $(q^2-1).$ Then
the other two exponents are either 0 or $(q^2-1),$ or else their
sum is $(q^2-1)$.   In these cases the $\lambda$-sums are all the
same for $\lambda'_0,\ldots,\lambda'_{2t-1}$. Using  (\ref{offorigin}), 
 the coefficient of
$T(xyz)^{q^2-1}$ becomes the constant term, which must be divisible by 
$q=p^\beta$, by Proposition~\ref{invariants} and the discussion following 
(\ref{offorigin}).
We are left 
with the cases that the $\Hh$-type 
of the monomial is either $(1,\ldots,1),$ or $ (2,\ldots,2)$.
If the $\Hh$-type is $(1,\ldots,1),$ there is nothing to prove.  We already showed 
(see (\ref{usum})) that the coefficients of monomials of $\Hh$-type
$(2,\ldots,2)$ are divisible by $q$.  

Now we consider terms in which $b'_{1},b'_{2},$ and $b'_{3}$ 
are all strictly between 0 and $q^2-1$.
We find it convenient at this point to go back to the affine
version of $\chi_{\overline U}$.  That is, we assume $x=1$ and do not 
write it.  We have
\begin{equation}
\nonumber
%\label{affchi}
\left(T(\alpha^qy^{2q})-T(\alpha y^2)+
T(\beta^q-\beta)T(y)^{q+1}-T(z)^q+T(z)\right)^{2(q^2-q)}.
\end{equation} 
We remember that each term has a nonzero power of $x$.
The only terms we have dropped are
$$\left(1-T(x)^{q^2-1}\right)T(y)^{q^2-1},$$ which we have 
already considered.
A typical term is 
\begin{equation}\label{term}
C\cdot D \cdot
T(y)^{2qk_1+2k_2+qk_3+k_3}T(z)^{qk_4+k_5},
\end{equation}
where $C=(-1)^{k_2+k_4}
T(\alpha)^{qk_1+k_2}T(\beta^q-\beta)^{k_3}$, 
$k_1+\cdots+k_5=2(q^2-q)$, and $D=\genfrac(){0 cm}0{2(q^2-q) }{ k_1,\ldots,k_5}$.
We will show that the multinomial coefficient is divisible by 
$p^\xi.$

Recall that the $p$-adic valuation of $D$
is the number of carries in
$k_1+\cdots+k_5=2(q^2-q)$.  So if
$$k_i=h_{i,2t}p^{2t}+\cdots+h_{i,1}p+h_{i,0}, \quad1\le
i\le5,$$ and $c_j$ represents the carry from the
$(j-1)^\mathrm{th}$ place to the $j^\mathrm{th}$ place, then
\begin{equation}\nonumber\begin{array}{rclclll}
h_{1,0}+\cdots+h_{5,0}&&-pc_1&=&0&&\\
h_{1,j}+\cdots+h_{5,j}&+c_{j}&-pc_{j+1}&=&0&\ \mathrm{for} \
&0<j<t\\ h_{1,\,t}+\cdots+h_{5,\,t}&+c_{t}&-pc_{t+1}&=&p-2&&\\
h_{1,j}+\cdots+h_{5,j}&+c_{j}&-pc_{j+1}&=&p-1&\ \mathrm{for}\
&t<j<2t\\ h_{1,2t}+\cdots+h_{5,2t}&+c_{2t}&&=&1,&&
\end{array}\end{equation} and $\nu_p (D)
=c_1+\cdots+c_{2t}.$ For some $j,\ 1\le j\le
t,$ assume that $c_j=c_{j+t}=0,$ so that
this position does not contribute to $\nu_p(D)$.
Clearly also $c_1=\cdots=c_{j-1}=0$ and we have 
$h_{i,0}=\cdots=h_{i,j-1}=0$ 
for each $i$.  Then
$$\sum_{\ell=0}^{j+t-1} p^{\ell} \sum_{i=1}^5 h_{i,\ell}=p^{t+j}-2q;$$
 $$(h_{1,j}+\cdots+h_{5,j})+
\cdots+(h_{1,j+t-1}+\cdots+h_{5,j+t-1})p^{t-1}=q-2p^{t-j};$$
 $$(h_{1,j+t}+\cdots+h_{5,j+t})
+\cdots+(h_{1,2t}+\cdots+h_{5,2t})p^{t-j}=2p^{t-j}-1.$$
Adding the last two expressions, multiplying by $(q+1)$, grouping the terms, and
using the formula for the twisted degrees of each $k_i,\ i=1,\ldots,5$, it 
immediately follows that
\begin{equation}\label{ksum}
\sum_{i=1}^5\Big((p^{t-j}k_i\mod q^2-1)
+(p^{2t-j}k_i\mod q^2-1)\Big)=q^2-1,
\end{equation} 
where each of the terms in (\ref{ksum}) is reduced before adding.

 We want 
to decide whether $s_j$ and $s_{j+t}$ can both be 2 in the $\Hh$-type 
of this monomial.  From (\ref{ssum}) we have, since 
$b'_{1}<q^2-1$,
\begin{eqnarray*}
E&=&(q^2-1)(s_j+s_{j+t}-2)\\
&=&
\sum_{i=1}^3
\left(
p^{2t-j}b'_{i} \mod (q^2-1)
+
p^{t-j}b'_{i} \mod (q^2-1)
\right)
\\&-&2(q^2-1)
\\&<&
\sum_{i=2}^3
p^{2t-j}b'_{i} \mod (q^2-1)
+
\sum_{i=2}^3 p^{t-j}b'_{i} \mod (q^2-1).
\end{eqnarray*}
Substituting for $b'_{2}$ and $b'_{3}$ from (\ref{term}), distributing the (mod 
$q^2-1$) operation, and using
(\ref{ksum}), we have
\begin{eqnarray*}
E&<&p^{2t-j}(2qk_1+2k_2+qk_3+k_3)\mod (q^2-1)
\\&+&p^{t-j}(2qk_1+2k_2+qk_3+k_3)\mod (q^2-1)
\\&+&
p^{2t-j}(qk_4+k_5)\mod (q^2-1)\\&+&
p^{t-j}(qk_4+k_5)\mod (q^2-1)\\
&\le&2\sum_{i=1}^3\Big((p^{t-j}k_i\mod q^2-1)
+(p^{2t-j}k_i\mod q^2-1)\Big)\\
&+&\sum_{i=4}^5\Big((p^{t-j}k_i\mod q^2-1)
+(p^{2t-j}k_i\mod q^2-1)\Big)\\
&<&2(q^2-1),
\end{eqnarray*}
where the last inequality is again strict because we assumed $b_3>0$.
We now have
$$s_j+s_{j+t}<4.$$
and the lemma is proved.  \qed

With Lemma~\ref{hval} and Lemma~\ref{bmval} we have proved the
 conjecture of Baker and Ebert in the Buekenhout-Metz case.

\begin{Theorem}
The number of points in the intersection of a Hermitian unital 
and a Buekenhout-Metz 
unital in $\PG(2,q^2)$ is congruent to 1 modulo $q$.
\end{Theorem}

Numerical evidence suggests that the Tits unital (in desarguesian planes of
order an odd power of two) also satisfy the conjecture of Baker and Ebert.  
In that case, Lemma~\ref{bmval} is not satisfied for individual terms of the 
expansion.

%%%%%%%%%%%%%%%%%%%%%%%%%%%%%%%%%%
\section{Two Hermitian varieties}\label{HvH}
Here we generalize Kestenband's result for two Hermitian unitals.
\begin{Theorem}
Let $H_1$ and $H_2$ be two nondegenerate Hermitian varieties in $\PG(n,q^2)$.  Then 
$$q^{n-1} \mid |H_1\cap H_2|.$$
\end{Theorem}
\Proof
Let the $\Hh$-types of $f$ and $g$ be $(s_0,\ldots,s_{2t-1})$ and
$(n+1-s_0,\ldots,n+1-s_{2t-1})$ and use (\ref{hermsum}).\qed

\end{document}